\documentclass{multi3}
\begin{document}
\chapter{Navier-Stokes Equations and Fluid Turbulence}
\author{P. Constantin}
\address{Department of Mathematics, The University of Chicago}
\begin{abstract}
An Eulerian-Lagrangian approach to incompressible
fluids that is convenient for both analysis and physics is presented.
Bounds on burning rates in combustion and heat transfer in 
convection are discussed, as well as results concerning spectra. 
\end{abstract}

\newcommand{\be}{\begin{equation}}
\newcommand{\ee}{\end{equation}}
\newtheorem{thm}{Theorem}

Incompressible fluids are described by the Navier-Stokes 
equation. Turbulence  (\cite{frisch}, \cite{sreeni}, 
\cite{tsinober}) experiments provide 
measurements that correspond to averages of certain 
quantities associated to 
the variables appearing in the Navier-Stokes equation. The 
present mathematical knowledge about the Navier-Stokes 
equations is incomplete. Some of the quantities measured 
in experiments are accessible to mathematical theory. They are 
usually low order, one-point bulk averages like the time 
average of integrals of squares of gradients. Most other 
measured quantities are not amenable to rigorous 
quantitative a priori analysis. Turbulence is concerned with 
statistical or collective properties of fluids. Nevertheless, 
the main impediment to progress in the rigorous analysis of 
turbulence is the present lack of understanding of possible 
blow up in individual solutions of the 
Euler and Navier-Stokes systems. 

I will discuss briefly the blow up problem
and present an Eulerian-Lagrangian approach to fluids.
I will also give examples of low order one-point bulk 
quantities that can be treated with present knowledge and
discuss results on certain two-point
quantities.

\section{An Eulerian-Lagrangian Approach to Fluids}

I will start by recalling that the Navier-Stokes-Euler system
can be written as an evolution equation for the three-component
velocity vector $u = u(x,t)$, 
$$
\frac{\partial u}{\partial t} + u\cdot\nabla u + \nabla p = 
\nu\Delta u + f;
$$
the pressure $p = p(x,t)$  preserves incompressibility
$$
\nabla \cdot u = 0.
$$
The Euler system is obtained if the 
kinematic viscosity vanishes, $\nu =0$; the Navier-Stokes system
if $\nu >0$. Boundary conditions are different for the two systems.
The blow up question can be stated in its simplest form
for the pure initial value problem: are there any smooth
initial data with finite energy that lead to solutions that diverge
in finite time? The answer is not known.

For a blow up in the Navier-Stokes system one would have to have
a finite time divergence of an 
eddy-viscosity-like quantity:
$$
\int_0^T\sup_{x, r}|u(x+r,t)-u(x,t)|^2dt <\infty\Rightarrow u\in{\mathbf
C}^{\infty}.
$$
By contrast, it is known (\cite{fgt})  that
$$
\int_0^T\sup_{x, r}|u(x+r,t)-u(x,t)|dt <\infty.
$$

In a situation in which all velocities are
finite no singularities can appear 
in the Navier-Stokes equations. 
One could accept the finiteness of velocities as a physical
hypothesis. For the Euler system this hypothesis would not be 
sufficient to ensure smoothness of solutions (\cite{maj}). A sufficient condition
for regularity in the Euler equations is the finiteness of the
time integral of the maximum magnitude of vorticity (\cite{bkm}). The
vorticity (the curl of velocity  or anti-symmetric part of
the velocity gradient), $\omega = \nabla\times u$ obeys a quadratic
evolution equation. The magnitude of vorticity obeys
$$
D_t|\omega| = \alpha |\omega|
$$
with $D_t = \partial_t + u\cdot\nabla$ the material derivative
along particle paths. The scalar stretching term $\alpha$ is related to 
the magnitude of vorticity by a principal value singular
integral (\cite{siam}):
$$
\alpha (x) = P.V. \int D\left (\hat{y},\xi (x),
\xi(x+y)\right )|\omega(x+y)|\frac{dy}{|y|^3}
$$
$$
\xi (x) = \frac{\omega(x)}{|\omega (x)|}.
$$
The smooth, mean zero function of three unit vectors 
$D$ vanishes when two of its arguments are
on the same line. Consequently, if the direction of vorticity
$\xi$ is Lipschitz then the singular integral representing
$\alpha$ is mild and the solutions remain smooth (\cite{ctm},
\cite{cfm}, \cite{cor}).
This is a generalization of the two dimensional situation where
$\nabla\xi = 0$ and the solutions remain smooth.
The possibility of blow up due to strong vortex stretching
is not removed by the previous result; the result only precludes
blow up in a smooth vortex line field. Blow up can occur also because of
strain intensification (the strain matrix is the symmetric part of the
gradient of vorticity). I will present now a description of
the Euler equations that is convenient for analysis (\cite{c-p}) and
allows for a clearer geometric picture of the possible singularity
formation. I will present the results only in the periodic case
for simplicity of exposition. The case of decay at infinity is
almost identical. 

\begin{thm} A function $u(x,t)$ solves the incompressible Euler 
equations if and only if it can be represented in the form 
$u = u_A$,
$$
u_A^i(x,t) = u_0^m(A(x,t))\frac{\partial A^m(x,t)}{\partial x_i}
 -\frac{\partial n_A(x,t)}{\partial x_i}
$$
where $A(x,t)$ 
solves the active vector equation
$$
\frac{\partial A(x,t)}{\partial t} + u_A(x,t)\cdot\nabla A(x,t) 
= 0
$$
with initial data
$$
A(x,0) = x.
$$
$u_0$ is the initial velocity 
and $n_A(x,t)$ is determined up to additive constants by the
requirement of incompressibility, $\nabla\cdot u_A = 0$.
\end{thm}

Let us denote by
$$
{\mathbf P} = \left ({\mathbf 1} - 
\nabla\Delta^{-1}\nabla\cdot\right)
$$
the Leray-Hodge projector on divergence-free vectors.
The local existence theorem requires just one derivative
to be H\"{o}lder continuous:

\begin{thm}
Let $u_0$ be a divergence free $C^{1,\epsilon}$ periodic vector
valued function of three variables. There exists a time
interval $[0,T]$ and a unique $C([0,T]; C^{1, \epsilon})$
spatially periodic vector valued function $\delta(x,t)$ such 
that
$$
A(x,t) = x + \delta (x,t)
$$
solves the active vector formulation of the Euler equations, 
$$
\frac{\partial A}{\partial t} + u_A\cdot\nabla A = 0,
$$
$$
u_A = {\mathbf {P}}\left \{(\nabla A(x,t))^*u_0(A(x,t))\right \}
$$
with initial datum $A(x,0) = x$.
\end{thm}

The proof of this result is based on an 
identity that removes the 
apparent ill-posedness, on singular integral calculus, and
on the use of the method of characteristics.
In the active vector formulation, the
"back-to-labels" map $A$ has conserved distribution.
Its time evolution is a smooth, volume-preserving
rearrangement. Singularities can occur only if the
gradient map $\nabla A$ diverges rapidly in finite time:
$$
\int_0^T\sup_x{|\nabla A(x,t)|^2dt} <\infty \Rightarrow
A\in{\mathbf C}^{\infty}.
$$

Thus, would-be singularities are gradient singularities in 
a conserved
quantity, similar to shocks in conservation laws, but with
the significant difference that the characteristic flow is measure
preserving. The formula relating the velocity to the value of the 
spatial gradient at the same instance of time (\cite{serr}, \cite{gold},
\cite{goldet}, \cite{hunt})  may have 
important mathematical consequences. On one hand, conservation of 
kinetic energy confers a constraint to the 
growth of $\nabla A$. On the other hand, the formula suggests 
that near regions of high gradient of $A$ the velocity is 
exceedingly high, making the regions of high gradient difficult 
to track and perhaps unstable. Dynamical stability or instability
of blow up modalities is a difficult subject. There are obvious
space and time symmetries (for instance, a minute delay of blow up), 
that clearly should not be categorized as instabilities. Nevertheless,
even relatively simple PDEs can exhibit the coexistence of
a variety of dynamical behaviors, including several
stable blow up modalities, stable time independent solutions, unstable
blow up modalities and dynamical connections between the unstable behaviors
and the stable ones (\cite{bckst}).

I will pass now from the blow up problem to some more tractable 
questions about average properties. One can obtain rigorous
upper bounds for certain bulk averages of solutions
of Navier-Stokes equations. Lower bounds are harder to obtain. 
Upper bounds for bulk averages of
low order moments for Rayleigh-B\'{e}nard convection will be 
described further below; I start with a lower bound for 
the burning rate in a simple model of turbulent combustion.

\section{Bulk Burning Rate}

Mixtures of reactants may interact in a burning region that 
has a rather complicated spatial structure but is thin across (\cite{Ro}).  
This reaction region moves towards the unburned reactants leaving 
behind the burned ones. When the reactants are carried by an 
ambient fluid then the burning rate is enhanced. The physical 
reason  for the observed speed-up is believed to be that 
fluid advection tends to increase the area available for reaction. 
What characteristics of the 
ambient fluid flow are responsible for burning rate enhancement? 
The question needs first to be made 
precise, because the reaction region may be 
complicated and, in general, may move with an ill-defined velocity.  
 An unambiguous quantity $V$  representing the bulk burning rate
is defined in (\cite{ckor}) and
explicit estimates of $V$ in terms of  
the magnitude of the advecting velocity and 
the geometry of streamlines are derived.  In situations where traveling waves are 
known to exist, $V$ coincides with the traveling wave speed and
the estimates thus provide automatically bounds 
for the speed of the traveling waves.  The main result of (\cite{ckor}) 
is the identification of a class of flows that are particularly  
effective  in speeding up the bulk burning rate. The main feature
of these  ``percolating flows''  is 
the presence of tubes of  streamlines connecting distant 
regions of burned and unburned material.  For such flows we obtained 
an optimal linear enhancement bound $ V\ge KU $ where $U$ represents 
the magnitude of the advecting velocity and $K$ is 
 a proportionality factor that depends on the geometry of
 streamlines but not the speed of the flow. Other flows and in
 particular cellular flows, which have  closed streamlines, on the 
other hand, may produce a weaker enhancement.  
The instantaneous bulk burning rate is defined by the formula
$$
V(t)=\int_D \frac{\partial T}{\partial t}(x,y,t) dxdy  
$$
where the integral extends over the spatial domain $D$, taken here for simplicity
of exposition to be a two-dimensional strip of unit width and infinite
length
$$
0\le y\le 1, \,\,-\infty <x< \infty.
$$
The temperature $T$ is assumed to obey Neumann boundary conditions
at the finite boundaries and to obey 
$$
T(-\infty, y)=1,
\,\,\,\,T(\infty, y)=0. 
$$ 
The simplified model is a passive reactive scalar with a
KPP nonlinearity
$$
 T_t + u \cdot \nabla T -\kappa \Delta T = \frac{v_{0}^{2}}{4\kappa}T(1-T). 
$$
with prescribed velocity $u$ that satisfies
$$ 
\int_0^1 u(x,y,t)\,dy = 0, \quad \nabla\cdot u = 0.
$$
The constant $v_0$ represents the speed of a stable one-dimensional
laminar ($u=0$) traveling wave. We start with a very general
lower bound:
\begin{thm} For arbitrary
initial data obeying
$$
0 \leq T_0(x,y) \leq 1.
$$
one has the general lower bound
\[ V(t) \geq Cv_0 \left( 1- e^{-\frac{v_{0}^{2}t}{2\kappa}} \right). \]
\end{thm}
The proof is based on the product lemma

\begin{lemma}
There exists a constant $C>0$ such that 
$$
0 \leq T(x,y) \leq 1,
$$
\[ T(-\infty, y)=1, \,\,\, T(\infty, y)=0 \,\,\, {\rm for \,\,\,any} \,\,\,y \in [0,1]. \]
implies
\[ \left( \int_D T(1-T) \, dxdy \right)
\left( \int _D|\nabla T|^2 \, dxdy \right) \geq C. \]
\end{lemma}

Although information about the velocity is not present in
the general result, it nevertheless shows that this model does not 
permit quenching. Also, the general lower bound applies to the 
homogenized version of the equations as well.

For a very general class of velocities $u(x,y,t)$ it can be shown that 
the bulk burning rate may not exceed a linear bound in  the 
amplitude of the advecting velocity.  For a large 
class of flows  we proved  lower bounds on the bulk burning rate that are 
 linear in the magnitude of advection. 
We denote by 
\[\langle V \rangle_\tau = \frac{1}{\tau} \int\limits_{0}^{\tau} V(t)
\,dt\]
the time average of the instantaneous bulk
burning rate.
The main result of (\cite{ckor}) is too technical to
state here precisely but its meaning is that
presence of coherent tubes of streamlines connecting 
unburned and  burned
regions enhances the burning rate 
$$
\langle V \rangle_{\tau}  \geq K U
$$                                            
as long as 
the velocity spatial scales 
are not too small compared to the reaction length scale
$\frac{\kappa}{v_0}$, and 
the time scale of change of the advecting velocity is not too small 
compared to  $\tau_0 = {\rm max}[\frac{\kappa}{v_0^2}, 
\frac{\tilde{H}}{v_0}]$ where $\tilde{H}$ is associated to
the width of the coherent tubes of streamlines.
For instance, a result concerning mean zero shear flow of the form
 \[ u(x,y)=(u(y),0), ~~\int_0^1 u(y)dy=0 \] can be stated 
as 
 \begin{thm}
  Let us consider an arbitrary partition of the interval $[0,1]$ into
  subintervals $I_j=[c_j-h_j,c_j+h_j]$ on which $u(y)$ does not change
  sign. Denote by $D_{-},$ $D_{+}$ the unions of intervals $I_j$ where $u(y)>0$
and $u(y)<0$ respectively. 
Then there exist  constants $C_{\pm} >0$, independent of the
  partition, $u(y),$ and the initial data $T_0(x,y)$, so that the
  average burning rate $\langle V \rangle_\tau$
satisfies the
  following estimate:
$$
  \langle V\rangle_\tau
 \ge C_+
c_+
\sum_{I_j \subset D_+}
\left(1+\frac{l^2}{h_j^2}\right)^{-1}
\int\limits_{c_j-\frac{h_j}{2}}^{c_j+\frac{h_j}{2}}|u(y)|dy 
$$
$$
+C_-c_-
\sum_{I_j \subset D_-}
\left(1+\frac{l^2}{h_j^2}\right)^{-1}
\int\limits_{c_j-\frac{h_j}{2}}^{c_j+\frac{h_j}{2}}|u(y)|dy
$$
for any $\tau \geq \tau_0=\hbox{max}\left[ \frac{\kappa}{v_0^2}, \frac{H}{v_0}
 \right].$ ($H=1$).
Here $l=\kappa/v_0.$ 
The constants $c_\pm$ are defined by 
\[ c_\pm = \left( \sum\limits_{I_j \subset D_\mp} \frac{h_j^3}{h_j^2 +l^2}
\right) \left( \sum\limits_{I_j} \frac{h_j^3}{h_j^2 +l^2}
\right)^{-1}. \]
\end{thm}

 The main result of (\cite{ckor}) applies to a large
class of flows that are not necessarily spatially periodic, nor shears, 
 and can  have completely arbitrary features outside the tubes of streamlines. 
 The bulk burning rate is still linear in the magnitude of
 the  advecting velocity, no matter what kind of behavior
 (closed streamlines, areas of still fluid, etc.) the flow has outside the
 tubes. The proportionality coefficient depends on the  geometry of  
the flow in a rather complex manner. These bounds can be 
extended to larger classes of chemistries. The lower bounds,
however, are not yet available for models in which 
there is a feedback coupling of temperature on the velocity. For such 
models upper bounds can be derived. In the next section we will discuss 
upper bounds for heat transfer. We will concentrate on the simplest
coupling, mediated by gravity in a Boussinesq approximation and
discuss the canonical case of Rayleigh-B\'{e}nard convection.

\section{Bulk Heat Transfer} 

The equations for Rayleigh-B\'{e}nard convection in the Boussinesq 
approximation are
$$
\frac{\partial u}{\partial t} + u\cdot\nabla u + \nabla p = \sigma
\Delta u + \sigma Ra \hat{e} T,
$$
$$
\nabla\cdot u = 0
$$
$$
\frac{\partial T}{\partial t} + u\cdot\nabla T = \Delta T.
$$
The vertical direction $\hat{e}$ of gravity is singled out.
We consider as spatial domain a box of height $1$ and lateral
side $L$. The velocity vanishes at the boundary, 
$T=1$ at the bottom boundary and $T = 0$ at top. The Nusselt number is
the space-time average of the flux of temperature across horizontal
cross-section planes. From the equations of motion it follows that 
$$
\left<|\nabla T|^2\right > = N.
$$
and also
$$
\left<|\nabla u|^2\right > = Ra(N
  - 1). 
$$
Here $<\cdots>$ is global space-time average. The general rigorous
result here is (\cite{con}, \cite{vection}) 
\begin{thm}
There exists an absolute constant $C$, independent of Rayleigh number
$Ra$, aspect ratio $L$ and Prandtl number $\sigma$ such that
$$
N \le 1 + C \sqrt{Ra}
$$
holds for all solutions of the Boussinesq equations.
\end{thm}

The experimental data (\cite{exp}) point however more to results of the type
$N\sim Ra^{\frac{1}{3}}$ or $N\sim  Ra^{\frac{2}{7}}$.
The exponent $\frac{1}{3}$ can be obtained rigorously for a
simplified model. Consider the infinite Prandtl number equations
for rotating convection 
$$
\left (\partial_t + u\cdot\nabla \right) T = \Delta T
$$
$$
\begin{array}{l} 
- \Delta u - E^{-1}v + p_x = 0 \\
-\Delta v  + E^{-1}u + p_y = 0 \\
-\Delta w  +  p_z    = RT.
\end{array}
$$
$$
\nabla\cdot {\mathbf {u}}= 0
$$
with boundary conditions: ($(u,v,w)$, $p$, $T$) 
periodic in  $x$ and $y$ with period $L$; $u$, $v$, and $w$ vanish 
for $z =0,1$, 
$T=0$ at $z=1$, $T=1$ at $z=0$. One can prove (\cite{chp})

\begin{thm} There exist absolute
constants $c_1,...,c_4$ so that the Nusselt number for rotating infinite 
Prandtl-number convection is bounded by 
$$
N - 1 \le 
$$
$$
\min \left \{c_1{R}^{\frac{2}{5}};\: (c_2E^2 + c_3E)R^{2};
 c_4 R^{1/3}(E^{-1} + \log_+ R)^{\frac{2}{3}}\right \}.
$$
\end{thm}

This coincides, in the limit of no rotation $E\to\infty$ with
 a logarithmic correction (\cite{ip}) to the  $\frac{1}{3}$ exponent. 
The bound also shows that strong rotation $E\to 0$ stabilizes the 
system and that increasing rotation may result in a 
non-monotonic behavior of the Nusselt number, as observed in experiments.

Consider now the horizontal average $\overline{T}(z,t)$ of
$T(x,y,z,t)$ and define
$$
n = \left <|\nabla (T-\overline{T})|^2\right >.
$$
Note that
$$
n \le N
$$

\begin{thm}
For the full Boussinesq system
$$
N\le 1 + c(n Ra )^{\frac{1}{3}}
$$
holds.
For the infinite Prandtl number system
$$
N\le 1 +C\left (Ra(\log_+Ra)^2\sqrt{n}\right )^{\frac{2}{7}}
$$
holds.
\end{thm}
Note that, if $n \le N$ is used then the result recovers the
exponents $\frac{1}{2}$ for general Rayleigh - B\'{e}nard
and logarithmically corrected $\frac{1}{3}$ for the infinite Prandtl
number case. But the rigorous appearance of the exponent $\frac{2}{7}$
is perhaps not coincidental.

One of the technical ingredients for the proof of these
results concerns zero order operators that are not translation 
invariant
$$
B = \frac{\partial^2}{\partial z^2}(\Delta^2_{DN})^{-1}\Delta_h
$$
where $w = (\Delta^2_{DN})^{-1}f$ is the solution of
$$
\Delta^2 w = f
$$
with horizontally periodic and vertically Dirichlet and Neumann
boundary conditions $w = w^{\prime} =0$. 

\begin{thm} For any $\alpha \in (0,1)$ there exists a positive
constant $C_{\alpha}$ such that every H\"{o}lder continuous
function $\theta$ that is horizontally periodic and
vanishes at the vertical boundaries satisfies
$$
\|B\theta\|_{L^{\infty}} \le C_{\alpha}\|\theta\|_{L^{\infty}}\left (1 +
\log_{+}\|\theta\|_{C^{0,\alpha}}\right )^2.
$$
\end{thm}
The proof of this result is based on
a pointwise bound on exponential-oscillatory sums 
of the type:
$$
K(x,z) = \sum_{k\in{\mathbf
    {Z}}^2}e^{\frac{2\pi}{L}ik\cdot{x}}m_k^pe^{-\epsilon m_k}
$$
where $\epsilon = \epsilon (z) \ge 0$,  
$m_k = \frac{2\pi}{L}|k|$ and 
$\epsilon (z) = 0\Rightarrow z=z_0$. The sum is singular at $z=z_0$
and the pointwise bound  
$$
|K(x,z)| \le C_p\left [|x|^2 +\epsilon^2(z)\right ]^{-\frac{p+2}{2}}
$$
is obtained using the Poisson summation formula.

The bounds on bulk one-point quantities presented above
are among the most successful areas of mathematical and experimental
agreement. The reason is perhaps that the quantities involved
are numbers, albeit numbers depending  on a parameter. The next step 
beyond the description of bulk one-point averages is the description
of power spectra. These are asserted to have some universal features
in physical turbulence theory; we present some mathematical results
in the next section.

\section{Spectra}

Unlike bulk one-point quantities, spectra are averages of functions. 
There are  some well-established spectra in the physical literature
associated to small scale turbulence: the Kraichnan spectrum in 
two dimensions and the Kolmogorov spectrum in three dimensions. 

The energy spectrum $E(k)$ is a function that has the property that
$$
\int_0^{\infty}E(k)dk = \left <|u|^2\right>.
$$
(The convention is that $<\cdots>$ is normalized 
space integral followed by long time average). The
Kolmogorov spectrum  for 3D turbulence is
$$
E(k) = C_{Kl}\left<\epsilon\right >^{\frac{2}{3}}k^{-\frac{5}{3}}.
$$
The 2D Kraichnan spectrum is 
$$
E(k) = C_{Kr}\left <\eta\right>^{\frac{2}{3}}k^{-3}.
$$
Here $\epsilon =\nu |\nabla u|^2$ is the rate of
dissipation of energy and $\eta = \nu |\nabla\omega |^2$
is the rate of dissipation of enstrophy. The spectra are supposed
to be  valid in a range of scales $k\in [k_i, k_d]$ where $k_d$
is the dissipation scale and is determined by viscosity and
$\epsilon$ (respectively viscosity and $\eta$) alone. Their
expressions are then determined by dimensional analysis. 
We consider the Littlewood-Paley decomposition of the velocity
associated to a mollifier $\phi (\xi)$. This is a smooth function
in $R^d$ that is  non-increasing, smooth, 
radially symmetric, satisfying $\phi (\xi) = 1$ for 
$|\xi |\le\frac{5}{8}$, $\phi (\xi) = 0$ for $|\xi| \ge\frac{3}{4}$. 
One sets $\psi_{(0)}(\xi)= \phi (\frac{\xi}{2}) - \phi(\xi)$, 
$\psi_{(m)}(\xi) = \psi_{(0)}(2^{-m}\xi)$ and 
$$
\phi (\xi) = \int_{{\mathbf R}^d}e^{-i(\xi\cdot z)}\Phi (z)dz,
$$
$$
\psi_{(m)}(\xi) = \int_{{\mathbf R}^d}
e^{-i(\xi\dot z)}\Psi_{(m)}(z)dz.
$$
The Littlewood-Paley decomposition is
$$
u(x,t) = u_{(-\infty )}(x,t) + 
\sum_{m = 0}^{\infty}u_{(m)}(x,t)
$$
where
$$
u_{(-\infty )}(x,t) = L^{-d}\int_{{\mathbf R}^d} 
\Phi \left (\frac{y}{L}\right )u(x-y,t)dy,
$$
$$
u_{(m)}(x,t) = L^{-d}\int_{{\mathbf R}^d}
\Psi_{(m)}\left (\frac{y}{L}\right ) u(x-y,t)dy
$$
and $L$ is a length (the integral scale).

We define the Littlewood-Paley spectrum to be
$$
E_{LP}(k) = k_m^{-1}\left <|u_{(m)}|^2\right >
$$
for $k_{m-1}\le k <k_{m}$, $m\ge 1$ with $k_m = 2^mL^{-1}$.

We start with the two dimensional Navier-Stokes equation
$$
\left (\partial_t + u\cdot\nabla -\nu\Delta\right )\omega = f
$$
with
$$
u (x,t) = \frac{1}{2\pi}\int_{{\mathbf R}^2}\frac{y^{\perp}}{|y|^2}\omega(x-y,t)dy
$$
One can prove (\cite{lit2})
\begin{thm} Assume that the source of vorticity in the Navier-Stokes
equations has spectrum localized in the  region of wave numbers
$k\le \frac{1}{L}$ for some $L>0$. Then there exists a constant $C$ such
that the Littlewood-Paley energy spectrum of solutions of two
dimensional forced Navier-Stokes equations obeys the bound
$$
E_{LP}(k) \le Ck^{-3}\left\{\tau^{-2}\left (\frac{k_d}{k}\right )^6 
\right\}
$$
for $k\ge L^{-1}$. Here $\tau^{-1} = <\|\nabla u\|_{L^{\infty}}>$.
\end{thm}

The corresponding three-dimensional energy spectrum result requires a
significant assumption:
$$
\langle |\nabla u|^3\rangle < \infty.
$$
Denoting
$$
\widehat{\epsilon} = \nu \left 
\{\langle |\nabla u|^3\rangle \right \}^{\frac{2}{3}}
$$
$$
\widehat{\eta} = \nu^{\frac{3}{4}}\left (\widehat{\epsilon}\right )^{-\frac{1}{4
}}
$$
$$
\widehat{k}_d = \nu^{-\frac{3}{4}}(\widehat{\epsilon})^{\frac{1}{4}}
$$
and setting
$$
C_{\psi} = \int\int |\nabla \Psi_{(0)} (a)||a|^2|\Psi_{(0)}(b)||b|dadb.
$$
we have (\cite{lit3}):
\begin{thm}
Consider three-dimensional body forces that satisfy
$$
\widehat{f}(k) = 0
$$
for all $|k|\ge \frac{C}{L}$ and some $C>0$.
Consider solutions of the three dimensional Navier-Stokes
equation that satisfy $\widehat{\epsilon} < \infty$.
Then
$$
E_{LP}(k) \le C_{\psi}\left (\widehat{\epsilon}\right )^{\frac{2}{3}}k^{-\frac{5
}{3}}
\left (\frac{k}{\widehat{k}_d}\right )^{-\frac{10}{3}}
$$
holds for $|k|\ge \frac{C}{L}$.
\end{thm}

{\bf Acknowledgment.} Research partially supported by NSF DMS-9802611, by AIM,
and by the ASCI Flash Center at the University of Chicago under DOE
contract B341495.

\bibliographystyle{plain}
 
\end{document}